\documentclass{article}
\usepackage[first]{draftcopy}
\usepackage{amsmath,amssymb,graphics}
\input xy
\xyoption {all}
\xyoption{v2}

\def\gg{\mathfrak {g}}
\def\Aoo{$A_\infty$}
\def\Loo{$L_\infty$}
\def\DHA{differential homological algebra}
\def\G{Gugenheim}

\newtheorem{definition}{Definition}

\newtheorem{remark}{Remark}

\newtheorem{theorem}{Theorem}

\begin{document}

\title{A Twisted Tale of Cochains and Connections}
\author{Jim Stasheff}

\maketitle

\large\centerline{In honor of  the 60-th birthday of Tornike Kadeishvilli}
\begin{abstract}
  Early in the history of higher homotopy algebra \cite{jds:hahII}, it was realized that Massey products are homotopy invariants in a special sense, but it was the work of Tornike Kadeisvili that showed they were but a shadow of an $A_\infty$-structure on the homology of a differential graded algebra. Here we relate his work to that of Victor Gugenheim \cite{victor:pert} and K.T. (Chester) Chen \cite{chen:iterated}. 
  This paper is a  personal tribute to Tornike and the Georgian school of homotopy theory as well as to Gugenheim and Chen, who unfortunately are not with us to appreciate this convergence.
\end{abstract}

\newpage
\tableofcontents
\section{Introduction}
Early in the history of higher homotopy algebra \cite{jds:hahII}, it was realized that Massey products are homotopy invariants in a special sense, but it was the work of Tornike Kadeisvili that showed they are but a shadow of an $A_\infty$-structure on the homology of a differential graded algebra. Here we relate his work to that of Victor Gugenheim \cite{victor:pert} and K.T. (Chester) Chen \cite{chen:iterated}.  However, in light of  the thorough technical analysis by  Huebschmann  elsewhere in this volume \cite{jh:tkfest} and his earlier survey in honor of Berikashvili \cite{jh:berikash}, this will be a more personal tribute to Tornike and the Georgian school of homotopy theory as well as to Gugenheim and Chen, who unfortunately are not with us to appreciate this convergence\footnote{Some of the memories and references here are my own, especially the personal ones, but I also owe a great deal to help from our Georgian colleagues and from Johannes Huebschmann, who also provided mathematical insights for earlier drafts.}.
 Essential to this discussion are the notions of  \emph{twisting element}  and special cases:
 \emph{twisting cochain} and \emph{flat connection}. In the latter guise, there are new applications on the border with mathematical physics with which we conclude this tribute
 to Tornike.
\section{Contact with the Georgian schools of category theory and algebraic topology}

Kadeishvili met Huebschmann at a meeting in Oberwohlfach in 1985, then he worked with Huebschmann  as a Humboldt Scholar  at Heidelberg in  1987/88.  (This collaboration produced, among other results, \cite{huebkade}.)  Meanwhile, Huebschmann, Gugenheim and I were at the International Topology Conference in Baku in 1987 and met several of the Georgian category and homotopy theorists, including Kadeishvili.
Even though we were able to meet in person in Baku, relations with the Soviet Union were still such that along the way to Baku I served as a courier between Borel at the IAS in  Princeton and Margolis in Moscow.

Since then, personal contact has increased in both directions with Kadeishvili and others visiting the West and major meetings  bringing non-Georgians to interact in Tiblisi.
Unfortunately, I have been unable to return myself, but two of my former students, Tom Lada and Ron Umble, have been my representatives. It was a pleasure to have Tornike in attendance at my 70th birthday (and Murray Gerstenhaber's 80th)  fest in Paris in 2007.

\section{ Berikashvili's twisting elements}
In 1968,  Berikashvili introduced the functor $\mathcal D$ 
in terms of \lq\lq twisting elements\rq\rq\   in a differential graded algebra $\mathcal A$. Such elements 
$\tau$ are  homogeneous and satisfy the 
equation
$$
d_{\mathcal A} \tau = \tau \tau,
$$
 where
$\tau \tau$
is the product in $\mathcal A$ of
$\tau$ with itself.

If, instead,  the algebra is a differential graded Lie algebra, the equation is
$$
d_{\mathcal A} \tau = 1/2[\tau, \tau].
$$
The element $\tau$ is necessarily of  degree $\pm 1$ (equal to that of $d_{\mathcal A} $) so
$[\tau, \tau]$ need not be zero.
Signs  depend on conventions and notation.

\begin{remark}
This equation has a long and honorable history in various guises. When the algebra is
that of differential forms on a Lie group,
it is called the \emph {Maurer-Cartan equation}. In deformation theory, it is the  \emph {integrability equation}. In mathematical physics, especially in the Batalin-Vilkovisky formalism, it is known as the  \emph {Master Equation}. 
At present, the name \emph {Maurer-Cartan equation} seems to have the upper hand. 
\end{remark}


Berikashvili's  functor
${\mathcal D}$ assigns to the dga $\mathcal A$ the set of equivalence classes of twisting elements, the equivalence relation nowadays being called \emph{gauge equivalence}. His point of view is quite appropriate to applications to deformation theory. On the other hand, his twisting elements include traditonal \emph{flat} connection forms
 $\omega$ on principal $G$-bundles: 
 $$
 G\to P \to M.
 $$
 That is, $\omega\in  \Omega^*(P,\gg),$
which is a differential graded Lie algebra for $\gg$, the Lie algebra of $G$.
 Modern generalizations include those initiated by Chen as well as those with $G$ generalized to \emph{higher structure} analogs. 
 
 \section{Connections}
\emph{Connections} have been well established for a long time in differential geometry, at least as far back as Elie Cartan, but, in the generality we build on,  the notion was introduced by Ehresmann \cite{ehresmann:connection},
though  at times the word refers to equivalent but distinct concepts. In particular, the word is sometimes used as a synonym for \emph{covariant derivative}.
\subsection{Ehresmann's connections}
\begin{definition}An \emph{Ehresmann connection} on a (locally trivial) smooth  fiber bundle
 $p:E\to B$ 
is a splitting of vector bundles for the induced morphism
$$TE \to B \times_B TB
$$
\noindent
of vector bundles. 
\end{definition}

\begin{remark}
In fact, Ehresmann's definition was for a submersion $p:E\to B$ and he proved that,
when the fibers are compact, such a splitting implies $p$ is actually locally trivial.
\end{remark}

To such a choice of horizontal subspaces in $TE$, there corresponds a \emph{connection form} 
$\omega$.
 Here we will be concerned mostly with a connection \emph{form} on something like a principal bundle. 
 Classically, for $\gg$ the  Lie algebra of a Lie group $G$
 and $\pi : P \to X$ a principal $G$-bundle, a \emph{principal
 connection form} on $P$ is a $\gg$-valued 1-form  $\omega \in
\Omega^1(P, \gg)$ which satisfies two conditions:
\begin{enumerate}
  \item
    $\omega$ restricts to the classical Maurer-Cartan $\gg$-valued 1-form on each
    fiber.
  \item
    $\omega$ is equivariant with respect to the adjoint $G$-action on $P$.
\end{enumerate}

Notice that a \emph{flat} principal Ehresmann connection form is a twisting element in the sense of Berikashvili.

\subsection{Cartan's connections}
Henri Cartan observed \cite{Cartan:g-alg} that  this could be expressed in terms of  a
morphism of graded-commutative algebras on which there is the action
of a Lie group (though only the action of the Lie algebra $\gg$ is
necessary). 

\begin{definition}A $\gg$\emph{-algebra} is a dgca (differential graded commutative algebra)  $A$ such that:

For each  $x\in \gg$, there is a  derivation called
`infinitesimal transformation' ${\mathcal L}(x)$ (today usually known
as the Lie derivative)  and a derivation called  `interior product' or `contraction'  $\imath (x)$
satisfying the relations:
\begin{enumerate}
  \item
   ${\mathcal L}\to Der A$ is an \emph{injective} dg Lie morphism 
  \item
  $\imath ([x,y]={\mathcal L}(x) \imath (y) - \imath (y){\mathcal L}(x) $
  \item
  ${\mathcal L}(x) = \imath (x)d + d\imath (x).$
\end{enumerate}

These derivations are respectively of degree 0 and degree $\pm 1$, opposite to the degree of $d$.
\end{definition}

These $\gg$-algebras are  also known as  \emph{Leibniz Pairs} \cite{fgv}.

The universal example of such a  $\gg$\emph{-algebra} is the Cartan-Chevalley-Eilenberg cochain algebra
$ \mathrm{CE}(\gg)$  for
Lie algebra cohomology:
$$ \mathrm{CE}(\gg) := \mathrm{Hom}(\Lambda \gg, \mathbf{R})$$
with the differential induced by extending the dual of the bracket as a derivation.

\begin{remark}
The originators expressed this in terms of alternating multilinear functions on $\gg$, which remains the correct formulation  for infinite dimensional Lie algebras, as opposed to the exterior algebra on the dual of
$\gg$.
\end{remark}

A \emph{Cartan connection}
$
  \xymatrix{
    \Omega^\bullet(P)
    &
    \mathrm{CE}(\gg)
    \ar[l]_{\omega}
  }
$
is then defined as respecting the operations $i(x)$ and ${\mathcal L}(x)$
for all $x\in\gg,$ but not necessarily respecting $d.$ 

If $\omega$ is a \emph{flat} connection,
it has curvature zero, that is equivalent to respecting $d$, hence satisfying the Maurer-Cartan equation.
Thus it is an example of a \emph{twisting element}.

 Once we are in the dg (differential graded) world, we could just as well take $\gg$ to be a differential graded Lie  algebra, using a completed tensor product $\hat\otimes$ where necessary. We can also work with differential graded associative algebras.  It was K.T. Chen who did this first in 1973 \cite{chen:iterated} and Kadeishvili independently in 1980 \cite{kad:Hfs}.
 
 \subsection{Chen's connections}
One of Chen's major contributions was a 
method for computing the real homology of the based loop space on a manifold in terms of the homology of 
the manifold. He effected this via his \emph{iterated integrals}, initially in \cite{chen:iterated} but evolving over several subsequent papers. In a very accessible survey \cite{chen:Bull77}, he uses the language of his \emph{formal power series connections}.


\begin{definition} Let $X$ be a graded vector space with basis $\{X_i\}$. A \emph{formal power series connection} on a differentiable space $M$ with values in a vector space $X$ is an element $\omega \in \Omega^*(M) [[X]]$ of the type 
$$
\omega  = \Sigma w_I X_I
$$
where $I$ denotes a multiindex $i_1\cdots i_r$ and  $X_I = X_{i_1}\cdots X_{i_r}$ and the coefficients 
$w_I$
are forms of positive degree on $M$.

\end{definition} 
The algebra $ \Omega^*(M) [[X]]$ can also be written as $ \Omega^*(M) \hat\otimes TX$ where
$TX$ is the tensor algebra on $X$.

Chen,
by suitably identifying his tensor product,
 saw that his
 condition for flatness becomes that of a twisting cochain, as he acknowledges in \cite{chen:Bull77}
Definition 3.2.1. 
In fact, such Chen connections with curvature zero are twisting elements in Berikashvili's sense,
though probably due to restricted communication with the Soviet Union, Chen did not reference Berikashvili.
 Contact between the western and USSR groups grew gradually, thanks  to the 
lifting of restrictions in the USSR under perestroika. Unfortunately, this came too late for Chen whose response to 
the ÒGeorgian schoolÓ we would very much like to have seen. 

To provide a multiplicative chain equivalence between his model and  the chains on the based loop space $\Omega X$, Chen  made use of his iterated integrals.  Thus his approach provided an `analytic' alternative to AdamsÕ cobar construction; one that was very useful in algebraic geometry
\cite{hain:ag}.



\section{Twisting cochains}
 The earliest occurance, to my knowledge, of the term  \emph{twisting cochain} is in the fundamental 1959 paper of E.H.Brown: Twisted tensor products I \cite{ehb:twist}. (In the S\'eminaire Henri Cartan 1956-57, there  is the term \emph{fonction tordante}, but that is in the context of simplicial sets, then known as `complete semi-simplicial complexes'.) Several related papers emphasized twisted tensor products and twisted differentials without mentioning twisting cochains, but it is the twisting cochains that are most closely related to connections. 
 
 \begin{definition} Given a  coaugmented differential graded coalgebra $C$ (with coaugmentation $\eta: R \to C$) and  an augmented differential graded 
algebra $A$ (with augmentation $\varepsilon\colon A \to R$, 
(both differentials being of degree $-1$),  a \emph{twisting cochain} $\tau: C \to A$ is a linear map of degree $-1$ satisfying the conditions 
$$
d_A \tau + \tau d_C = \tau\smile\tau
$$
$$
 \varepsilon \tau=0 \ \ \rm{and}\ \   \tau \eta = 0.$$
The Òcup-productÓ $\smile$ is defined in the module of linear maps $ C \to A$  by using the coproduct 
 $\Delta$ on $C$ and the product $m$ on $A$: Given two maps $ f , g : C \to A$, 
$$f \smile g = m(f \otimes g)\Delta$$. 
 \end{definition}

  Again, if we take $\mathrm{Hom}(C,A)$ with cup product as $\mathcal A$, then twisting cochains are twisting elements in Berikashvili's sense. 
  
  By 1960, Gugenheim had access to a preprint of Brown's paper and became interested in the idea of a twisting cochain and its relation to the description of a simpicial fibre bundle as a `twisted cartesian product' in his work with Barratt and Moore \cite{bgm}.

The fundamental role of twisting cochains in \emph{differential homological algebra} was developed by J.C. Moore \cite{jcm:DHA}. In 1974, Moore, together with Husemoller and Stasheff, emphasized this role and applied it to a `classical' problem in algebraic topology \cite{hms}. Readers of that paper may well surmise who had primary responsibility for which part.

 
 \subsection{Chen's Theorem}
 In 1973, Chen \cite{chen:loop} proved a result which can be paraphrased as follows.  Let  $\Omega M$ denote the based loop space on $M$ and $T(s^{-1}H_\bullet(M))$ denote the tensor algebra on the \emph{desuspension} of the vector space $H_\bullet(M)$.
 
 \begin{theorem}
 For a simply connected manifold $M$, there is a twisting element $\omega \in \Omega^\bullet(M)\hat\otimes T(s^{-1}H_\bullet(M))$ with respect to  a derivation $\partial$ on $T(s^{-1}H_\bullet(M))$ for which there is a map 
 $$\Theta: C_\bullet(\Omega M) \to (T(s^{-1}H_\bullet(M)), \partial)$$ giving an isomorphism in homology.
 \end{theorem}
 \begin{remark} $ \Omega^\bullet(M)\hat\otimes T(s^{-1}H_\bullet(M))$ can be written as 
 $ \Omega^\bullet(M, T(s^{-1}H_\bullet(M)).$
 \end{remark}
 
Soon after, \G \ focused on the fact  that twisting cochains and homotopies of twisting cochains are at the heart of Chen's 
work. This interest culminated in 1982 \cite{victor:pert} where Gugenheim gave an algebraic version of Chen's 
theorem on the homology of the loop space, not restricted to the smooth setting and differential forms nor even real coefficients. For this, $C_\bullet( M)$  is replaced by a suitable differential graded coalgebra $C$ and $\Omega C$ denotes Adams' cobar construction on $C$. \G\ constructed a multiplicative perturbation $\partial$ of the cobar differential 
on $\Omega H(C)$ and a map 
$\Omega C\to \Omega_\partial H(C)$ which is a purely algebraic 
analog of the map $\Theta$ given by Chen's iterated integrals. 
Also in the early 80s, Huebschmann made extensive use of twisting cochains and homological perturbations, \cite{jh:origins} and references therein. For the early history of homological perturbation theory (HPT), the review {\tt{MR1103672}} of \cite{GLS:II}
 by Ronnie Brown is excellent.

There were independent developments in the USSR by Berikashvili, Kadeishvili, Saneblidze and 
others.
\subsection{Kadeishvili's theorem}
In 1980 and quite independently,  Kadeishvili proved the corresponding very basic result
for algebras and the bar construction, which is denoted  $B$ and  generalized to apply to \Aoo-algebras where needed:

\begin{theorem}
If $A$ is an augmented differential graded algebra with $H_{\bullet}(A)$ free as a module over the ground ring, 
then $H_{\bullet} (A)$ admits an \Aoo-structure such that there exists a map of dgcoalgebras 
$BH (A)\to BA$ 
inducing an isomorphism in homology. 
\end{theorem}

This result is sometimes referred to as a `minimality theorem', which, I think,  has the wrong emphasis and point of view. It is the transfer of structure up to homotopy that to me is most important.

Apparently \Aoo-structures caught on faster in Moscow and especially Tiblisi than in the US, where
\G's version  came to Stasheff's attention. 
In 1986 \cite{gugjds} together they made the connection with \Aoo-structures.
Considerable ÒwesternÓ work  was thus inspired by Chen's ideas, whereas Berikashvili and  Kadeishvili led the way in the ``east''.

Once \Aoo-structures appeared in this context, it was natural to consider $A$ itself being an \Aoo-algebra; this is what Kadeishvili did in 1982 \cite{tk:1982}. He also developed further the relation between \Aoo-structures and Massey products, which was only  implict in my early work.

One of the characteristic features of Chen's connections and Gugenheim's twisting cochains is that they 
include as special cases twisted tensor products which are acyclic. 


\section{The Lie and \Loo versions and mathematical physics}

We can also consider a twisting function $\tau:C\to L$ from a dg coalgebra to a dg Lie algebra. 
As far as I know, this first occurred in the context of rational homotopy theory in Quillen's seminal paper
\cite{quillen:qht}.

The main advantage of using $\mathrm Hom$ is the manifest naturality and the avoidance of finiteness conditions. Similarly the originators of Lie algebra cohomology got it right: using alternating multilinear functions on a Lie algebra $\gg$ rather than the exterior (Grassmann) algebra on the dual of $\gg$.
This works for infinite dimensional Lie algebras as well. Connection forms with values in a Lie algebra
play key roles in math and physics, so generalizations to values in an \Loo-algebra are natural.


 Just as a Lie algebra or dg Lie algebra $\gg$ can be characterized by a `quadratic' differential on the graded symmetric coalgebra on the (de)suspension of $\gg$, so \Loo-algebras can be characterized by removing  the quadratic restriction.
However, there was a considerable lag in introducing and developing \Loo-structures until they were needed in algebraic deformation theory \cite{Sjds} and string field theory \cite{z:csft}. They were however implicit  in Sullivan's models in rational homotopy theory \cite{sullivan:inf} in 1977. A lot was going on that year!



Chen deals with Lie algebras primarily in his studies of fundamental groups, but Hain \cite{hain:tc}
adapts Chen's twisting cochains in the form of twisting elements in $A\hat\otimes L$
where $L$ is a dg Lie algebra (with $d$ of degree -1) and $A$ is a dg commutative algebra 
(with $d$ of degree +1). As a student of Chen's and with specific computational examples in mind,
this setting is natural for Hain. Generalization to \Loo-valued connections is needed for applications to mathematical physics. 

\subsection{Bundles with \Loo-structure}

 As higher category theory was developed, mimicing homotopy theory, Lie 2-algebras (also known as infinitesimal crossed modules)  appeared \cite{baez-crans:Lie2alg} and were recognized as special (very small) \Loo-algebras.  This led naturally to the \DHA\  version of classical differential geometry, in particular, generalized connections, curvature and `all that'. However, the driving force in this recent development was application suggested by mathematical physics: differential graded string theory and even `5brane theory' \cite{SSS:5brane}.
 
 The first example of  ``higher bundles with connection''
occurred with the fundamental (super)string coupling to the Neveu-Schwarz
(NS) $B$-field  $B_2  $\cite{Neveu-Schwarz}.
This $B$-field  is a connection on a 2-bundle
and appears in an action functional
$\int_{\Sigma} B_2$ for the
string worldsheet (surface) $\Sigma$.
Here a 2-bundle \cite{bartels:thesis} means a bundle with fibres which are at least 2-vector spaces \cite{baez-crans:Lie2alg}, that is,
a differential graded vector space of the form $V=V_0 \oplus V_1$ with $d$ of degree 1.
Similarly, 
a connection $B_6$ on a 6-bundle
appears in an action functional $\int_{\Sigma_6} B_6$
for the fivebrane worldvolume $\Sigma_6$.
This and related matters are  explained in \cite{Freed} in the language of differential characters and in
\cite{SSS:Loo, SSS:5brane} in the language of higher bundles.


To stay at the level of dgcas, we make the definition below. First notice that the definition of a $\gg$-algebra above applies to \Loo-algebras $\gg$ except that, for $x$ of degree $k$, the degree of 
$\mathfrak{L}(x)$ is $ - k$ and that of $i(x)$ is $-k-1$.

\begin{definition}An \emph{algebra of differential forms on a principal  $\gg$-\emph{bundle} over a smooth space $X$} is a $\gg$-algebra in the sense of H. Cartan (we denote it $\Omega^\bullet(Y)$) with a monomorphism $\pi: \Omega^\bullet(X)\to \Omega^\bullet(Y)$ such that $i(x)\pi = 0 = \mathfrak{L}(x)\pi$
 for all $x\in\gg$.
 \end{definition}

The formal definition of such a connection can be stated in the Henri Cartan form. The graded commutative algebra $\mathrm{CE}(\gg)$ has the usual operations $i(x)$ and ${\mathcal L}(x)$.
As before, an \emph{ \Loo-Cartan connection}
$
  \xymatrix{
    \Omega^\bullet(Y)
    &
    \mathrm{CE}(\gg)
    \ar[l]_{A}
  }
$
is then defined as a graded algebra map injective on  the dual of $\gg$, respecting  $i(x)$ and ${\mathcal L}(x)$
for all $x\in\gg,$ but not necessarily respecting $d.$



In most of the applications to physics, $\gg$ is non-zero in a very, very small number of degrees.
These are related to connected covers of $BO$ or $BU.$ For example, a \emph{spin}-structure 
on a smooth space $X$ corresponds to a lifting of the classifying map of $TX$ to the 2-connected
cover, a \emph{string}-structure corresponds to a lifting to the 4-connected cover and  a \emph{Fivebrane}-structure corresponds to a lifting to the 8-connected cover. The standard Chern-Weil approach using the \Loo\  version of the Weil algebra then applies to determine characteristic classes of bundles with such structures \cite{SSS:Loo}.

\subsection{Higher spin structures, closed string field theory and \Loo-algebra}

Recognition that the mathematical structure of  sh-Lie algebras (=  \Loo-algebra) was
appearing in physics first occurred in my discussions at UNC with Burgers (visiting van Dam)
and then
with Zwiebach at the third GUT Workshop in 1982.
In their study of \emph{field dependent gauge symmetries} for field theories for higher spin particles 
\cite{burgers:diss,BBvD:probs,BBvD:arb}, Behrends, Burgers and van Dam discovered what turned out to be an \Loo-structure.
In conversations, Burgers and I found we had common formulas, if not a common language. 


As a generalization of Lie algebras, sh Lie algebras (now more commonly known as \Loo-algebras)
appeared in physics as symmetries or gauge transformations, though they were not presented as such initially in the physics literature \cite{burgers:diss,BBvD:arb,z:csft}.  They were recognized as such in closed string field theory when Zwiebach and I were together at the third GUT Workshop in Chapel Hill.
The corresponding Lagrangians consist of (sums of)
$(N+1)$\emph{-point functions}
They can be regarded as being formed from the $N$-fold brackets
$ \lbrack x_1,x_2,\dots ,x_N\rbrack $ of the \Loo-algebra
by evaluation with a dual field via an inner product. 
In terms of the $N$-fold bracket, we then define
$$
\{ y_0 y_1\dots y_N\} = <y_0\vert \lbrack y_1,y_2,\dots ,y_N\rbrack > .
$$
\par
Zwiebach presents a
classical action in closed string
field theory,  gauge transformations
and shows the invariance of the action.
The classical string action is simply given by
$$S(\Psi ) = {1\over 2} \langle \Psi , Q\Psi \rangle
+ \sum_{n=3}^\infty {\kappa^{n-2} \over n!}
\{ \Psi\dots\Psi \}.
$$
The gauge transformations of the theory are given by
$$\delta_\Lambda \vert \Psi >
= \sum_{n=0}^\infty { \kappa^n \over n!}
 [\Psi,\dots,\Psi, \Lambda ] $$
 for $\Lambda$ in $\gg.$
Notice that all the terms of higher order are necessary for these to be
consistent.
\par
Similarly, in open string field theory,
one can define $N+1$-point functions
using the structure maps $m_N$ of an \Aoo-algebra.  

A primer on \Loo theory for physicists is \cite{ls}.

\subsection{Open-closed string field theory and OCHA}
Having considered both \Aoo- and \Loo-algebras, we come to the combination known as \emph{OCHA} for Open-Closed Homotopy Algebra \cite{kaj-jds:math,kaj-jds:physics}. Inspired by open-closed string field theories \cite{zwiebach:mixed}, these involve an \Loo-algebra acting by derivations (up to strong homotopy) on an \Aoo-algebra but have an additional piece of structure corresponding to a closed string opening to an open string. The details are quite complicated in the original papers, but, just as other ``$\infty$'' algebras can be characterized by a single coderivation on an appropriate dgc coalgebra, the same has been achieved for OCHAs by Hoefel \cite{hoefel}. Now, returning to Kadeishvili's work, in a recent paper
 \cite{kad-lada}, he and Lada have exhibited a very small, concrete example, providing one that perhaps can provide a \emph{toy model} for open-closed string field theory.


\section{Coda}
It has been a pleasure to sketch the connections between Tornike's work and my own, as well as that of many other contributers to higher homotopy algebra. Surely some further `twist' to this history lies ahead.




\section{References}
\begin{thebibliography}{BBvD86}

\bibitem[AN71]{Neveu-Schwarz}
J.H.~Schwarz A.~Neveu.
\newblock Factorizable dual model of pions.
\newblock {\em Nucl.Phys.B}, 31:86--112, 1971.

\bibitem[Bar04]{bartels:thesis}
Toby Bartels.
\newblock Higher gauge theory {I}: 2-bundles.
\newblock {\tt{arXiv:math/0410328}}, 2004.

\bibitem[BBvD85]{BBvD:probs}
F.A. Berends, G.J.H. Burgers, and H.~van Dam.
\newblock On the theoretical problems in constructing interactions involving
  higher spin massless particles.
\newblock {\em Nucl.Phys.B}, 260:295--322, 1985.

\bibitem[BBvD86]{BBvD:arb}
F.A. Berends, G.J.H. Burgers, and H.~van Dam.
\newblock Explicit construction of conserved currents for massless fields of
  arbitrary spin.
\newblock {\em Nucl.Phys.B}, 271:429--441, 1986.

\bibitem[BC04]{baez-crans:Lie2alg}
John~C. Baez and Alissa~S. Crans.
\newblock Higher-dimensional algebra. {VI}. {L}ie 2-algebras.
\newblock {\em Theory Appl. Categ.}, 12:492--538 (electronic), 2004.

\bibitem[BGM59]{bgm}
M.~G. Barratt, V.~K. A.~M. Gugenheim, and J.~C. Moore.
\newblock On semisimplicial fibre-bundles.
\newblock {\em Amer. J. Math.}, 81:639--657, 1959.

\bibitem[Bro59]{ehb:twist}
Edgar~H. Brown, Jr.
\newblock Twisted tensor products. {I}.
\newblock {\em Ann. of Math. (2)}, 69:223--246, 1959.

\bibitem[Bur85]{burgers:diss}
G.J.H. Burgers.
\newblock {\em On the construction of field theories for higher spin massless
  particles}.
\newblock PhD thesis, Rijksuniversiteit te Leiden, 1985.

\bibitem[Car50]{Cartan:g-alg}
H.~Cartan.
\newblock Notions d'alg\'ebre diff\' erentielle; application aux groupes de
  {L}ie et aux vari\'et\'es o\u' op\`ere un groupe de {L}ie.
\newblock In {\em Colloque de Topologie, Bruxelles (1950)}, pages 15--27. CBRM,
  1950.

\bibitem[Che73a]{chen:iterated}
K.-T. Chen.
\newblock Iterated integrals of differential forms and loop space homology.
\newblock {\em Ann.~of~Mathematics}, 97:217--246, 1973.

\bibitem[Che73b]{chen:loop}
K.-T. Chen.
\newblock Iterated integrals of differential forms and loop space homology.
\newblock {\em Ann.~of~Mathematics}, 97:217--246, 1973.

\bibitem[Che77]{chen:Bull77}
Kuo~Tsai Chen.
\newblock Iterated path integrals.
\newblock {\em Bull. Amer. Math. Soc.}, 83(5):831--879, 1977.

\bibitem[Ehr51]{ehresmann:connection}
Charles Ehresmann.
\newblock Les connexions infinit\'esimales dans un espace fibr\'e
  diff\'erentiable.
\newblock In {\em Colloque de topologie (espaces fibr\'es), {B}ruxelles, 1950},
  pages 29--55. Georges Thone, Li\`ege, 1951.

\bibitem[FGV95]{fgv}
M.~Flato, M.~Gerstenhabe, and A.~A. Voronov.
\newblock Cohomology and deformation of {L}eibniz pairs.
\newblock {\em Lett. Math. Phys.}, 34:77--90, 1995.

\bibitem[Fre00]{Freed}
Daniel~S. Freed.
\newblock Dirac charge quantization and generalized differential cohomology.
\newblock In {\em Surveys in differential geometry}, Surv. Differ. Geom., VII,
  pages 129--194. Int. Press, Somerville, MA, 2000.
\newblock {\tt{hep-th/0011220v2}}.

\bibitem[GLS91]{GLS:II}
V.~K. A.~M. Gugenheim, L.~A. Lambe, and J.~D. Stasheff.
\newblock Perturbation theory in differential homological algebra. {II}.
\newblock {\em Illinois J. Math.}, 35(3):357--373, 1991.

\bibitem[GS86]{gugjds}
V.K.A.M. Gugenheim and J.D. Stasheff.
\newblock On perturbations and ${A}_{\infty}$-structures.
\newblock {\em Bull. Soc. Math. Belg.}, 38:237, 1986.

\bibitem[Gug82]{victor:pert}
V.K.A.M. Gugenheim.
\newblock On a perturbation theory for the homology of a loop space.
\newblock {\em J. Pure and Appl. Alg.}, 25:197--207, 1982.

\bibitem[Hai83]{hain:tc}
Richard~M. Hain.
\newblock Twisting cochains and duality between minimal algebras and minimal
  {L}ie algebras.
\newblock {\em Trans. Amer. Math. Soc.}, 277(1):397--411, 1983.

\bibitem[Hai02]{hain:ag}
Richard Hain.
\newblock Iterated integrals and algebraic cycles: examples and prospects.
\newblock In {\em Contemporary trends in algebraic geometry and algebraic
  topology ({T}ianjin, 2000)}, volume~5 of {\em Nankai Tracts Math.}, pages
  55--118. World Sci. Publ., River Edge, NJ, 2002.

\bibitem[HK91]{huebkade}
Johannes Huebschmann and Tornike Kadeishvili.
\newblock Small models for chain algebras.
\newblock {\em Math. Z.}, 207(2):245--280, 1991.

\bibitem[HMS74]{hms}
Dale Husemoller, John~C. Moore, and James Stasheff.
\newblock Differential homological algebra and homogeneous spaces.
\newblock {\em J. Pure Appl. Algebra}, 5:113--185, 1974.

\bibitem[Hoe06]{hoefel}
Eduado Hoefel.
\newblock On the coalgebra description of{ OCHA}.
\newblock {\tt{arXiv:math/0607435}}, 2006.

\bibitem[Hue99]{jh:berikash}
J.~Huebschmann.
\newblock Berikashvili's functor $\mathcal{D}$ and the deformation equation\/.
\newblock {\em Festschrift in honor of N.~Berikashvili's 70th birthday,
  Proceedings of A.~Razmadze Institute}, 119:59--72, 1999.
\newblock {\tt{math.AT/9906032}}.

\bibitem[Hue09a]{jh:tkfest}
J.~Huebschmann.
\newblock On the construction of ${A}_\infty$-structures.
\newblock {\em Festschrift in honor of T.~Kadeishvili's 60th birthday}, 2009.
\newblock {\tt{arXiv:0809.4791}}.

\bibitem[Hue09b]{jh:origins}
J.~Huebschmann.
\newblock Origins and breadth of the theory of higher homotopies\/.
\newblock In {\em Festschrift in honor of M. Gerstenhaber's 80-th and Jim
  Stasheff's 70-th birthday}, Progress in Math. (to appear). 2009.
\newblock {\tt{arxiv:0710.2645}}.

\bibitem[Kad80]{kad:Hfs}
T.V. Kadeishvili.
\newblock On the homology theory of fibre spaces.
\newblock {\em Russian Math. Surv.}, 35:3:231--238, 1980.
\newblock {\tt{math.AT/0504437}}.

\bibitem[Kad82]{tk:1982}
T.~V. Kadeishvili.
\newblock The algebraic structure in the homology of an ${A}(\infty )$-algebra.
\newblock {\em Soobshch. Akad. Nauk Gruzin. SSR}, 108:249--252, 1982.

\bibitem[KL]{kad-lada}
T.V. Kadeishvili and T.~Lada.
\newblock A small open-closed homotopy algebra ({OCHA}).
\newblock {\em Georgian Mathematical Journal}.
\newblock to appear.

\bibitem[KS06a]{kaj-jds:math}
Hiroshige Kajiura and Jim Stasheff.
\newblock Homotopy algebras inspired by classical open-closed string field
  theory.
\newblock {\em Comm. Math. Phys.}, 263(3):553--581, 2006.

\bibitem[KS06b]{kaj-jds:physics}
Hiroshige Kajiura and Jim Stasheff.
\newblock Open-closed homotopy algebra in mathematical physics.
\newblock {\em J. Math. Phys.}, 47(2):023506, 28, 2006.

\bibitem[LS93]{ls}
T.~Lada and J.D. Stasheff.
\newblock Introduction to sh {L}ie algebras for physicists.
\newblock {\em Intern'l J. Theor. Phys.}, 32:1087--1103, 1993.

\bibitem[Moo71]{jcm:DHA}
John~C. Moore.
\newblock Differential homological algebra.
\newblock In {\em Actes du Congr\`es International des Math\'ematiciens (Nice,
  1970), Tome 1}, pages 335--339. Gauthier-Villars, Paris, 1971.

\bibitem[Qui69]{quillen:qht}
D.~Quillen.
\newblock Rational homotopy theory.
\newblock {\em Annals of Mathematics}, 90:205--295, 1969.

\bibitem[SS77]{Sjds}
M.~Schlessinger and J.~Stasheff.
\newblock Rational homotopy theory -- obstructions and deformations.
\newblock In {\em Proc. Conf. on Algebraic Topology, Vancouver}, pages 7--31,
  1977.
\newblock LMM 673.

\bibitem[SSS09a]{SSS:5brane}
Hisham Sati, Urs Schreiber, and Jim Stasheff.
\newblock Fivebrane structures.
\newblock {\tt{ arXiv:0805.0564}}, 2009.

\bibitem[SSS09b]{SSS:Loo}
Hisham Sati, Urs Schreiber, and Jim Stasheff.
\newblock ${L}_\infty$-algebra connections and applications to {S}tring- and
  {C}hern-{S}imons n-transport.
\newblock {\em TBA}, 2009.
\newblock {\tt{ arXiv:0801.3480}}.

\bibitem[Sta63]{jds:hahII}
J.~Stasheff.
\newblock Homotopy associativity of {H}-spaces, {II}.
\newblock {\em Trans. Amer. Math. Soc.}, 108:313--327, 1963.

\bibitem[Sul77]{sullivan:inf}
D.~Sullivan.
\newblock Infinitesimal computations in topology.
\newblock {\em Pub. Math. IHES}, 47:269--331, 1977.

\bibitem[Zwi93]{z:csft}
B.~Zwiebach.
\newblock Closed string field theory: {Q}uantum action and the
  {B}atalin-{V}ilkovisky master equa\-tion.
\newblock {\em Nucl. Phys. B}, 390:33--152, 1993.

\bibitem[Zwi98]{zwiebach:mixed}
B.~Zwiebach.
\newblock Oriented open-closed string theory revisited.
\newblock {\em Ann. Phys.}, 267:193--248, 1998.
\newblock {\tt{hep-th/9705241}}.

\end{thebibliography}
\end{document}